\documentclass[11pt]{article}
\usepackage{amsfonts,amsmath,amssymb,amsthm, amscd, subcaption}
\usepackage{graphicx}
\usepackage[stable]{footmisc}
\numberwithin{equation}{section}

\usepackage{tikz}
\usepackage{tikz-cd}
\usetikzlibrary{decorations.markings}
\usetikzlibrary{arrows}
\usetikzlibrary{calc}

\newtheorem{theorem}{Theorem}[section]

\newtheorem{lemma}[theorem]{Lemma}

\theoremstyle{definition}

\usepackage{float}

\def\({\langle \hskip -.1cm \langle} 
\def\){\rangle \hskip -.1cm \rangle} 
\def\<{{\langle}}
\def\>{{\rangle}}

\def\S{{\mathbb S}}

\newcommand{\Mr}{M_A^{\textup{red}}}

\allowdisplaybreaks

\def\ni{\noindent} 

\begin{document}

\title{Peripheral elements in reduced Alexander modules: an addendum} 

\author{Daniel S. Silver
\and
Lorenzo Traldi}

\maketitle 


\begin{abstract} We answer a question raised in ``Peripheral elements in reduced Alexander modules'' [J. Knot Theory Ramifications {\bf 31} (2022), 2250058]. We also correct a minor error in that paper.
\end{abstract}

\begin{center} Keywords: \textit {checkerboard coloring, classical link, longitude, Seifert surface, virtual link};   AMS Classification: MSC: 57K10 \end{center}

\section{Introduction} \label{Intro} 

In \cite[Defs.\ 1 and 2]{T}, Traldi introduced a modified version of the traditional peripheral structure for the reduced Alexander module $\Mr(L)$ of a classical or virtual link $L=K_1 \cup \dots \cup K_\mu$, with each link component $K_i$ having a set $M_i(L)$ of associated meridians and a single associated longitude $\chi_i(L)$. The module $\Mr(L)$ is generated by $\cup M_i(L)$, and the longitudes of $L$ lie in the \emph{reduced Alexander invariant} of $L$, the submodule of $\Mr(L)$ consisting of sums $\sum \lambda_j m_j$ where the $m_j$ are meridians and $\sum \lambda_j=0$. The longitude of a classical or virtual knot is 0, and the set of meridians is determined up to an automorphism of the reduced Alexander invariant, so the peripheral structure is of little value when $\mu=1$. For links of two or more components, though, the peripheral structure significantly increases the sensitivity of the reduced Alexander module as a link invariant.

The main purpose of this addendum to \cite{T} is to prove that the equality 
\[
\sum_{i=1}^ \mu \chi_i(L)=0
\]
holds for all classical links, answering a question posed in \cite{T}. The proof is a simple modification of a Seifert surface argument that seems to be the standard way to prove the well-known fact that the longitude of a classical knot lies in the second commutator subgroup of the fundamental group of the complement. (This fact is equivalent to the assertion that the longitude is $0$ in $\Mr(L)$.) See \cite[p.\ 39]{BZ03}, for instance.

To keep the addendum brief we do not repeat definitions of terms that appear in \cite{T}, or are standard in the literature.

Given an oriented classical link $L$, it is well known that the infinite cyclic cover (also called the ``total linking number cover'') $\widetilde X$ of the exterior $X=\S^3-L$ may be constructed using a Seifert surface $S$ as follows. First, cut the $3$-sphere $\S^3$ along $S$ to obtain a cobordism $W$ between copies $S_-, S_+$ of $S$. Then, glue countably infinitely many copies of $W$ end-to-end, by identifying $S_+$ in the $i$th copy of $W$ with $S_-$ in the $(i+1)$st copy. See \cite[p.\ 54]{BZ03} or \cite[p.\ 129] {rolfsen} for a detailed discussion of the construction when $L$ is a knot; there is no significant difference when $L$ has several components. Notice that if we choose a basepoint $*$ for $X$ in the interior of $S$, and use the preimage of $*$ on the $0$th copy of $S_+$ as a basepoint for $\widetilde X$, then $S$ lifts to a homeomorphic copy $\widetilde S$ in $\widetilde X$.

It is also well known that if $p:\widetilde X \to X $ is the covering map then the reduced Alexander module $\Mr(L)$ is isomorphic to the relative homology group $H_1(\widetilde X, p^{-1}(*))$, and the reduced Alexander invariant corresponds to the image  $i_*(H_1(\widetilde X)) \subset H_1(\widetilde X, p^{-1}(*))$, where $i:\widetilde X \to (\widetilde X, p^{-1}(*))$ is the inclusion map. See for instance Crowell's account in \cite{C3}.

Recall that in \cite{T}, a longitude is chosen for each link component $K_i$, and the total linking number with $L$ of each of these longitudes is $0$. Then each longitude of $L$ represents one component of the oriented boundary of $\widetilde S$ in $H_1(\widetilde X)$. The sum of these longitudes is $0$ in the homology group $H_1(\widetilde X)$ for the simple reason that the sum represents a cycle bounded by $\widetilde S$.

\section{Virtual links}

It is easy to find examples of virtual links for which the sum of longitudes is not $0$ in $\Mr$. One such example is the virtual Hopf link, discussed in \cite[Sec.\ 6]{T}. It has only one nonzero longitude.  

Suppose for the moment that a virtual link $L$ can be imbedded in a thickened surface $\Sigma \times [0,1]$ in such a way that $L$ has an oriented spanning surface $S$. That is, $L$ is \emph{almost classical} in the sense of Silver-Williams \cite{SW} and Boden {\sl et al.\ }\cite{BCK}. Then the argument of Section 1 may be applied to the coned-off exterior $X=(\Sigma \times [0,1] -L )/\Sigma \times 1$ and an infinite cyclic cover $\widetilde X$ of $X$. The fact that the reduced Alexander invariant coincides with the first homology group $H_1(\widetilde X)$ can again be verified as in \cite{C3}.

If $L$ has an unoriented spanning surface $S$ -- that is, $L$ is {\sl checkerboard colorable} -- then $S$ cannot be used in the same way to construct an infinite cyclic cover of $X$. However, $S$ can be used in an analogous way to construct a $2$-fold cover $X_2$, with $\pi_1(X_2)$ the subgroup of $\pi_1(X)$  consisting of elements represented by words in the meridional generators with even total exponent. Suitably modified, the argument of Section 1 then implies that the sum of longitudes is $0$ in the quotient module $\Mr(L)/(t+1)\Mr(L)$. 

\section{Two additional comments on \cite{T}}

There is an algebraic error in the statement of \cite[Lemma 8]{T}. The corrected statement, given below, is proven by a simple calculation using the definitions from \cite{T}.  The corrected statement has the same logical significance as the flawed statement given in \cite[Lemma 8]{T}: the peripheral structure of a classical or virtual link $L$ determines all the linking numbers in $L$. 
\begin{lemma}
\label{linklem}
Let $L=K_1 \cup \dots \cup K_{\mu}$ be a classical or virtual link. If $i \in \{1, \dots, \mu\}$, then $\epsilon^\mu \varphi_L (\chi_i(L))$ is the element of $\mathbb Z ^ \mu$ whose $j$th coordinate is $0$ if $j=1$, $\ell_{j/i}(K_i,K_j)$ if $i \neq j \neq 1$, and 
\[
-\sum _{k \neq i} \ell_{k/i}(K_i,K_k)
\]
if $i=j \neq 1$.
\end{lemma}

We should also mention that the classical case of \cite[Theorem 9]{T}, which states that the longitudes in $\Mr(L)$ generate the submodule annihilated by $t-1$, is equivalent to a result of Sakuma \cite[Theorem 4]{S}.

\ni Department of Mathematics and Statistics\\
\ni University of South Alabama\\ Mobile, AL 36688 USA\\
\ni Email: silver@southalabama.edu
\bigskip

\ni Department of Mathematics\\
\ni Lafayette College\\Easton, PA 18042 USA\\
\ni Email: traldil@lafayette.edu\\ 

\end{document}